\documentclass{article}

\usepackage[ansinew]{inputenc}
\usepackage{amsmath}
\usepackage{amssymb}
\usepackage{amsthm}
\usepackage{amscd}
\usepackage{amsfonts}
\usepackage[english]{babel}
\usepackage{graphicx}
\newtheorem {defn}{Definition}

\newtheorem {prop}[defn]{Proposition}
\newtheorem {cor}[defn]{Corollary}
\newtheorem {thm}[defn]{Theorem}

\newtheorem {lemma}[defn]{Lemma}

\newcommand{\tr}{\operatorname{tr}}
\newcommand{\diag}{\operatorname{diag}}
\newcommand{\bbmatrix}[1]{\left[ \begin{array}{cccc} #1 \end{array} \right]}
\newcommand{\SP}{{\operatorname{sp}}}
\newcommand{\SPEC}{{\operatorname{spec}}}
\newcommand{\IKdV}{\boldsymbol{\mathcal{I}}_{KdV}}
\newcommand{\A}{\boldsymbol{\mathcal{A}}}

\title{Solution of the KdV equation using evolutionary vessels}
\author{Andrey Melnikov\\Drexel University}

\begin{document}
\maketitle
\abstract{In this work we present a new method for solving of the Korteweg-de Vries (KdV)
equation
\[ q'_t = - \dfrac{3}{2} q q'_x + \dfrac{1}{4} q'''_{xxx}.
\]
The proposed method is a particular case of the theory of evolutionary vessels, developed
in this work.
Inverse scattering of the Sturm-Liouville operator and evolution of its potential are
the basic ingredients, similar to the existing methods developed
by Gardner-Greene-Kruskal-Miura (1967), Zacharov-Shabbath (1974) and
Peter Lax (1977). Evolutionary KdV vessel may be considered as a generalization of 
these previous works. The advantage of the new method is that it produces a unified approach
to existing solutions of the KdV equation. For example, odd or even 
analytic, periodic, almost periodic solutions are shown to be particular cases of this theory.
Generalizing this method we can also produce many
PDEs, associated with integrable systems, in an arbitrary number of variables
(in the spirit of Zakarov-Shabat).
}

\section{Introduction}
The Korteweg-de Vries (KdV) is the following nonlinear evolutionary PDE
\begin{equation} \label{eq:KdV}
q'_t = - \dfrac{3}{2} q q'_x + \dfrac{1}{4} q'''_{xxx}.
\end{equation}
The standard technique for solving this equation \cite{bib:GGKM}
consists of the detailed study of the evolution of the scattering data
corresponding to the potential $q(x,t)$ of Sturm-Liouville differential equation
\begin{equation}\label{eq:SL}
 -y''_{xx} + q(x,t) y = \lambda y.
\end{equation}
Zakarov-Shabat \cite{bib:ZakSha74} have created a similar scheme for a wide family
of nonlinear equations using the ideas of P. Lax
\cite{bib:LaxPair}, or more precisely a Lax pair.
We will present an alternative scheme, similar to these ideas, for constructing
solutions of the KdV equation \eqref{eq:KdV}. 

The starting point is a $2\times 2$ matrix-function of a
complex variable $\lambda$, realized \cite{bib:bgr} in the following form using an auxiliary Hilbert
space $\mathcal H$:
\[ \begin{array}{lll}
S(\lambda) = I - B_0^* \mathbb X_0^{-1} (\lambda I - A)^{-1} B_0 \sigma_1, \\
A \mathbb X_0 + \mathbb X_0 + A^* + B_0 \sigma_1 B_0^* = 0, \quad \mathbb X_0^* = \mathbb X_0,
\end{array} \]
where $\sigma_1 = \bbmatrix{0&1\\1&0}$ and $B_0:\mathbb C^2\rightarrow\mathcal H$,
$A, \mathbb X_0: \mathcal H\rightarrow\mathcal H$ are linear operators. These kinds of realizations were
studied for example in \cite{bib:bgr, bib:BL, bib:MyThesis, bib:SLVessels, bib:GenVessel, bib:Potapov, 
bib:SahACSpec, bib:SchurPont, bib:Brodskii}. The main feature in the study of such realizations is the
correspondence between the multiplicative structure of the matrix $S(\lambda)$ and the invariant subspaces
of the operator $A$ \cite{bib:BL, bib:Potapov}.
The pioneering work \cite{bib:defVess} of M. Livsi\v c
and its generalization in \cite{bib:Vortices} are the main sources of motivation for this research project.

The variable $\lambda$ may be called the
``frequency'' variable. Define next
$\sigma_2 = \bbmatrix{1&0\\0&0}$, $\gamma=\bbmatrix{0&0\\0&i}$ and solve the following
differential equations with initial conditions in the order they appear, starting from $B_0, \mathbb X_0$
appearing in the realization of $S(\lambda)$ above:
\[ \begin{array}{lllll}
0  =  \frac{d}{dx} (B(x)\sigma_1) + A B(x) \sigma_2 + B(x) \gamma & B(0) = B_0, & \text{\eqref{eq:DB}} \\
\frac{d}{dx} \mathbb X(x)  =  B(x) \sigma_2 B^*(x), & \mathbb X(0) = \mathbb X_0, & \text{\eqref{eq:DX}} \\
 \dfrac{\partial}{\partial t} B(x,t) = i A \dfrac{\partial}{\partial x} B(x,t),
 & B(x,0) = B(x), & \text{\eqref{eq:DBt}} \\ 
\dfrac{\partial}{\partial t} \mathbb X(x,t) = 
 i A B(x,t) \sigma_2 B^*(x,t) - i B(x,t)\sigma_2 B^*(x,t) A^* +
 i B(x,t)\gamma B^*(x,t), & \mathbb X(x,0) =\mathbb X(x). & \text{\eqref{eq:DXt}} \\
\end{array} \]
Then we can check that the Lyapunov equation
\[ A \, \mathbb X(x,t)+ \mathbb X(x,t) A^* + B(x,t) \sigma_1 B(x,t) = 0
\]
holds for all $x,t$ (the derivative with respect to $x$ and $t$ is zero). The so-called ``tau function''
of this construction is
\[ \tau(x,t) = \det(\mathbb X_0^{-1} \mathbb X(x,t))
\]
and has the property that $q(x,t) = -2 \dfrac{\partial^2}{\partial x^2}\tau(x,t)$ satisfies 
the KdV equation \eqref{eq:KdV} (Theorem \ref{thm:KdVequation}). 
In Theorem \ref{thm:Solitons} it is shown that when the function
$S(\lambda)$ is meromorphic (i.e. has a finite-dimensional realization, $\mathcal H<\infty$), we obtain a soliton.
The collection of these constructed operators and spaces 
\[
\mathfrak{E}_{KdV} = (A, B(x,t), \mathbb X(x,t); \sigma_1, 
\sigma_2, \gamma, \gamma_*(x,t);
\mathcal H, \mathbb C^2).
\]
is called a KdV evolutionary vessel (see Definition \ref{def:KdvVessel}).

In this work we show that the scheme proposed here for constructing solutions of the
KdV equation is effective and that some important solutions such as solitons, periodic, almost
periodic, and analytic are obtained in this way. Moreover, for the case of
odd/even in $x,t$ solutions of \eqref{eq:KdV} we show in Section \ref{sec:Examples} that
there is also a 1-1 correspondence with some ``canonical'' constructions of the KdV vessels
$\mathfrak{E}_{KdV}$.

Section \ref{sec:SL} is devoted to the vessel, obtained using only the first two steps
\eqref{eq:DB} and \eqref{eq:DX}, and as a result it deals with operators depending on $x$ only.
These types of vessels were studied by the author and others in a series of papers
\cite{bib:MyThesis,bib:SLVessels,bib:SchurVessels,bib:GenVessel} and they represent a basis for this theory.
The definition and basic ideas of vessels follow from the work \cite{bib:defVess} of M. Livsi\v c.
In Section \ref{sec:Scatt} it is shown that the classical inverse scattering fits into this
setting and that the remaining two equations \eqref{eq:DBt}, \eqref{eq:DXt} are equivalent to
the evolutions of the potential $q(x,t)$, or of its scattering data as in the classical case.

In Section \ref{sec:Examples} we prove that there is 1-1 correspondence between certain odd at $t=0$
classes of solutions and vessels.
These correspondences are presented in Table \ref{table:KdV}.
\begin{table}[htc]
\hspace{1cm}\begin{tabular}{|c||c|c|c|}
\hline
Class of solutions & Type of evol. vessel & Previous works & Correspondence \\
\hline
	Solitons & $\dim\mathcal H<\infty$ & \cite{bib:BKSWSolitons} & Theorem \ref{thm:Solitons} \\
	\hline 
	Analytic-   & Regular vessels & \cite{bib:BDT} & Theorem \ref{thm:AnalytBetaReal} \\
	exponential &&& \\
	\hline
	$\exists q'''_{xxx}(x,t)$ on $\mathbb R$ & $\SPEC(A)\subseteq i\mathbb R_+$ & \cite{bib:BKSWSolitons, bib:LaxPair}, & Theorem  \ref{thm:L1Sol} \\
	&  & \cite{bib:ZakSha74, bib:GGKM} & \\
	\hline 
	$T$-periodic in $x$ & $\mathcal H=\ell^2$, $A=\diag(\sqrt{-1}n_i^2)$, & \cite{bib:LaxPeriod},
	& Theorem \ref{thm:PeriodSol}\\
	$\dfrac{T^3}{(2\pi)^2}$-periodic in $t$ &  $n_i=\dfrac{2\pi N_i}{T}, N_i\in\mathbb N$  &  
	\cite{bib:Krich77} & \\
	\hline 
	almost periodic & $\mathcal H=\ell^2$, $A=\diag(\sqrt{-1}r_i^2)$, & \cite{bib:LaxAlmPeriod} & Theorem \ref{thm:AlmostPeriodSol} \\
									& $r_i\in\mathbb R$ & & \\	
\hline
\end{tabular}
\caption{Classes of KdV odd solutions and corresponding vessels.}
\label{table:KdV}
\end{table}

It seems, and is a project by itself, that the solution of G. Segal, G. Wilson using
Loop groups \cite{bib:SegalWilson} corresponds to some choices of the initial function
$S(\lambda)$. The Krichever solutions \cite{bib:Krich77}, which can also be implemented using 
loop groups \cite{bib:SegalWilson}, are conjecturally examples of vessels 
implemented on curves \cite{bib:GenVessel}. Baker functions \cite{bib:Baker28} seem to be
generalized by the transfer function of the vessel
\[ S(\lambda,x,t) = I - B^*(x,t)\mathbb X^{-1}(x,t)(\lambda I - A)^{-1} B(x,t)\sigma_1.
\]
Finally, the Zakarov-Shabat scheme \cite{bib:ZakSha74} can be generalized by adding 
variables to the operator $B,\mathbb X$ and defining linear differential equations in the spirit
of \eqref{eq:DB}, \eqref{eq:DX}, \eqref{eq:DBt}, and \eqref{eq:DXt}.
For example, we should be able to produce solutions of the Kadomtzev-Petviashvili equation \cite{bib:KadPetv}
using this scheme.

\section{\label{sec:SL}Sturm-Liouville vessels}
\subsection{Sturm-Liouville Vessels}
\begin{defn} \label{def:SLparam}
The Sturm Liouville (SL) vessel parameters are defined as follows
\[ \sigma_1 = \bbmatrix{0 & 1 \\ 1 & 0},
\sigma_2 = \bbmatrix{1 & 0 \\ 0 & 0},
\gamma = \bbmatrix{0 & 0 \\ 0 & i},
\gamma_*(x)=\bbmatrix{-i (\beta'(x) - \beta^2(x)) & -\beta(x) \\ \beta(x) & i}
\]
for a real-valued differentiable function $\beta(x)$, defined on an interval $\mathrm I$.
\end{defn}
Before we define a notion of a vessel which involves an auxiliary Hilbert space $\mathcal H$ and operators
(for $x\in\mathrm I$)
\begin{equation} \label{eq:DefOperats} \begin{array}{llllllll}
A, \mathbb X(x) &:& \mathcal H &\rightarrow & \mathcal H, \\
B(x)						&:& \mathbb C^2 &\rightarrow & \mathcal H
\end{array} \end{equation}
we need to consider some regularity assumptions.
We assume that the operator $A$ may be unbounded with a domain $D(A)$.
Moreover, certain algebraic and differential relations connect these operators. 
As as a result, we must determine assumptions, which ensure
that the relations between $A, \mathbb X(x), B(x)$ become solvable equations.
\begin{defn}[Regularity assumptions] \label{def:UnOperts}
Operators $A, \mathbb X(x), B(x)$ are said to satisfy regularity assumptions on $\mathrm I$ if
\begin{enumerate}
	\item $A$ is a generator of $C_0$ semi-group,
	\item $B(x) \mathbb C^2 \in D(A^n)$ for all $x\in\mathrm I,n\in\mathbb N$,
	\item The operator $\mathbb X(x)$ is self-adjoint and invertible for all $x\in\mathrm I$.
\end{enumerate}
\end{defn}

\begin{defn} 
A Strum-Liouville vessel is a collection of operators and spaces
\begin{equation} \label{eq:DefSLV}
\mathfrak{V} = (A, B(x), \mathbb X(x); \sigma_1, 
\sigma_2, \gamma, \gamma_*(x);
\mathcal{H},\mathbb C^2), 
\end{equation}
The space $\mathcal H$ is Hilbert and the 
operators $A, \mathbb X(x), B(x)$ are defined in \eqref{eq:DefOperats} so
that the regularity assumptions hold. The operators are subject to the following \textbf{vessel conditions}:
\begin{align}
\label{eq:DB} 0  =  \frac{\partial}{\partial x} (B\sigma_1) + A B \sigma_2 + B \gamma, \\
\label{eq:XLyapunov} A \mathbb X + \mathbb X A^*  =   B \sigma_1 B^*, \\
\label{eq:DX} \frac{\partial}{\partial x} \mathbb X  =  B \sigma_2 B^*, \\
\label{eq:Linkage}
\gamma_*  =  \gamma + \sigma_2 B^* \mathbb X^{-1} B \sigma_1 -
 \sigma_1 B^* \mathbb X^{-1} B \sigma_2.
\end{align}
\end{defn}
By definition, the \textbf{transfer function} of this vessel is
\[ 
S(\lambda,x) = I - B^*(x) \mathbb X^{-1}(x) (\lambda I - A)^{-1} B(x) \sigma_1.
\]
Poles and singularities of $S$ with respect to $\lambda$ are determined by $A$ only.
Multiplication by the function $S(\lambda,x)$ maps \cite{bib:MyThesis, bib:GenVessel, bib:MelVin1} 
solutions $u(\lambda, x)$ of the input LDE with the spectral parameter $\lambda$
\begin{equation}
\label{eq:InCC}
 \lambda \sigma_2 u(\lambda, x) -
\sigma_1 \frac{\partial}{\partial x}u(\lambda,x) +
\gamma u(\lambda,x) = 0
\end{equation}
to solutions $y(\lambda, x) = S(\lambda,x) u(\lambda, x)$ of the output LDE with the same spectral parameter
\begin{equation} \label{eq:OutCC}
\lambda \sigma_2 y(\lambda, x) - \sigma_1 \frac{\partial}{\partial x}y(\lambda,x) +
\gamma_*(x) y(\lambda,x) = 0.
\end{equation}

Denoting $u_\lambda(x) = \bbmatrix{u_1(\lambda,x)\\ u_2(\lambda,x)}$ we shall obtain that the input compatibility
condition \eqref{eq:InCC} is equivalent to
\[ \left\{ \begin{array}{lll}
-\frac{\partial^2}{\partial x^2} u_1(\lambda,x) = -i\lambda u_1(\lambda,x), \\
u_2(\lambda,x) = - i \frac{\partial}{\partial x} u_1(\lambda,x).
\end{array}\right.
\]
The output $y_\lambda(x) = \bbmatrix{y_1(\lambda,x)\\y_2(\lambda,x)} = S(\lambda,x) u_\lambda(x)$ satisfies the output
equation \eqref{eq:OutCC}, which is equivalent to
\[ \left\{ \begin{array}{lll}
-\frac{\partial^2}{\partial x^2} y_1(\lambda,x) + 2 \beta'(x) y_1(\lambda,x) = -i\lambda y_1(\lambda,x), \\
y_2(\lambda,x) = - i [ \frac{\partial}{\partial x} y_1(\lambda,x) - \beta(x) y_1(\lambda,x)].
\end{array}\right.
\]
In other words, we can see that multiplication by $S(\lambda,x)$ maps the solution of the trivial SL equation 
\eqref{eq:SL} ($q=0$)
to solutions of the more complicated one, defined by the potential 
\begin{equation} \label{eq:qbeta}
q(x) = 2 \beta'(x).
\end{equation}
As a result, the function $S(\lambda,x)$ satisfies the following differential equation
\begin{equation} \label{eq:DS}
\dfrac{\partial}{\partial x} S = 
\sigma_1^{-1}(\lambda \sigma_2 + \gamma_*) S - S \sigma_1^{-1}(\lambda \sigma_2 + \gamma).
\end{equation}

We notice that $(\mathcal H, \mathbb X(x))$ form a Krein space: $\mathcal H$
is a set, whose (Krein) inner product $[u,v]=\langle\mathbb X(x)u,v\rangle, \forall u,v\in\mathcal H$ depends on $x$ and is differentiable. 
From system theory \cite{bib:Brodskii, bib:Kalman} and operator theory related to $J$-contractive functions
\cite{bib:bgr, bib:Potapov} we adopt some of the following additional characterizations of the vessel.
\begin{defn} \label{defn:ClassVess}
The vessel $\mathfrak{V}$ \eqref{eq:DefSLV} is called
\begin{itemize}
	\item \textbf{dissipative}, if $\mathbb X(x)>0$ for all values of $x\in\mathrm I$,
	\item \textbf{Pontryagin}, if $\mathbb X(x)$ has $\kappa\in\mathbb N$ negative squares at the
		right half plane for all values of $x\in\mathrm I$,
	\item \textbf{Krein}, if $\mathbb X(x)$ has infinite number of negative squares at the
		right half plane for all values of $x\in\mathrm I$,
	\item \textbf{regular},  if all the operators $A, B(x), \mathbb X(x)$ are bounded operators for all $x$,
	\item \textbf{minimal}, if for all $x$,
		\begin{equation} \label{eq:minCond}
		 \overline{\operatorname{span}} \{ A^n B(x) \mathcal E\mid n\in\mathbb N \} = \mathcal H.
		\end{equation}
\end{itemize}
\end{defn}

\subsection{\label{sec:StandConstr}Standard construction of a vessel}
in this section we show that one can easily construct vessels. Choose two Hilbert spaces
$\mathcal H, \mathcal E$ and define three
operators $\mathbb X_0, A:\mathcal H \rightarrow \mathcal H$ and $B_0:\mathcal E\rightarrow\mathcal H$
such that $\mathbb X_0$ is invertible and the following equalities hold:
\[ \mathbb X_0^* = \mathbb X_0, \quad
A \mathbb X_0 + \mathbb X_0 A^* +  B_0 \sigma_1 B_0^* = 0.
\]
Then solve \eqref{eq:DB} with the initial value
\[ 0  =  \frac{d}{dx} (B(x)\sigma_1) + A B(x)  \sigma_2 + B(x) \gamma, 
\quad B(x_0) =  B_0,
\]
and solve equation \eqref{eq:DX}:
\[ \mathbb X(x)  =  \mathbb X_0 + \int_{x_0}^ x B(y) \sigma_2 B(y)^* dy.
\]
Finally, define $\gamma_*(x)$ from $\gamma(x)$ using \eqref{eq:Linkage}. Thus a vessel is created \cite{bib:GenVessel}:
\begin{lemma} The collection
\[ \mathfrak{K_V} = (A, B(x), \mathbb X(x); \sigma_1, 
\sigma_2, \gamma, \gamma_*(x);
\mathcal{H}, \mathcal{E};\mathrm I)
\]
is a vessel.
\end{lemma}
\noindent\textbf{Proof:} Follows from the definitions. \qed

We can obtain in this manner a rich family of vessels, using standard models, which create operators
$\mathbb X_0, A, B_0$:
\begin{enumerate}
	\item Liv\v sic model of a non-selfadjoint operator \cite{bib:BL}, where $\mathbb X=I$, $A+A^*+B J B^* = 0$, and
		$J$ is a signature matrix ($J=J^*=J^{-1}$),
	\item Theory of nodes, developed in \cite{bib:Brodskii},
	\item Krein space realizations for symmetric functions \cite{bib:KreinReal},
	\item Vessels on curves \cite{bib:GenVessel},
\end{enumerate}
The most important aspect of all these constructions is that the singularities of the transfer
function are manageable. Under the minimality condition (see Definition \ref{defn:ClassVess}), taking the operator $A$ as a
diagonal operator, we obtain that these singularities are at certain points 
($\dim\mathcal H<\infty$) or on a countable set of points ($\mathcal H=\ell^2$)
or on a curve $\Gamma$ ($\mathcal H = L^2(\Gamma)$). In this manner we can construct 
vessels whose transfer function will have singularities at the prescribed points (or
on prescribed curves). The converse of this construction is an analogue of inverse scattering.

\subsection{The tau function of a vessel}
Following the ideas presented in \cite{bib:SLVessels}, we define the tau function of the vessel
$\mathfrak V$, (see \eqref{eq:DefSLV}) as follows
\begin{defn} \label{def:Tau} For a given vessel
\[
\mathfrak{V} = (A, B(x), \mathbb X(x); \sigma_1(x), 
\sigma_2(x), \gamma(x), \gamma_*(x);
\mathcal{H}, \mathcal{E};\mathrm I),
\]
the tau function $\tau(x)$ is defined as
\begin{equation} \label{eq:Tau} \tau = \det (\mathbb X^{-1}(x_0) \mathbb X(x))
\end{equation}
for an arbitrary point $x_0\in\mathrm I$.
\end{defn}
Notice that, using vessel condition \eqref{eq:DX}, $\mathbb X(x)$ has the formula
\[ \mathbb X(x) = \mathbb X(x_0) + \int\limits_{x_0}^x B^*(y) \sigma_2 B(y) dy,
\]
and as a result
\[ \mathbb X^{-1}(x_0) \mathbb X(x) = I + \mathbb X^{-1}(x_0) \int\limits_{x_0}^x B^*(y) \sigma_2 B(y) dy.
\]
Since $\sigma_2$ has rank 1, 
this expression is of the form $I + T$, for a trace class operator $T$. Since 
$\mathbb X_0$ is an invertible operator, there exists a non-trivial interval (of length at least $\dfrac{1}{\|\mathbb X_0^{-1}\|}$) on which $\mathbb X(x)$ and $\tau(x)$ are defined. Recall \cite{bib:GKintro} that a function $F(x)$ from (a, b) into the group G (the set of bounded invertible operators on H of the form $F(x)=I + T(x)$) is said to be 
differentiable if $F(x) -I$ is \textit{differentiable} as a map into the trace-class operators. In our case,
\[ \dfrac{d}{dx} (\mathbb X^{-1}(x_0)\mathbb X(x)) = 
\mathbb X^{-1}(x_0) \dfrac{d}{dx} \mathbb X(x) =
\mathbb X^{-1}(x_0) B(x)\sigma_2B^*(x)
\]
exists in trace-class norm. Israel Gohberg and Mark Krein \cite[formula 1.14 on p. 163]{bib:GKintro}
proved that if $\mathbb X^{-1}(x_0)\mathbb X(x)$ is a differentiable function
into G, then $\tau(x) = \SP (\mathbb X^{-1}(x_0)\mathbb X(x))$
\footnote{$\SP$ - stands for the trace in the infinite dimensional space.} is a differentiable map into $\mathbb C^*$ with
\begin{multline} \label{eq:GKform}
\dfrac{\tau'}{\tau}  = \SP (\big(\mathbb X^{-1}(x_0)\mathbb X(x)\big)^{-1} 
\dfrac{d}{dx} \big(\mathbb X^{-1}(x_0)\mathbb X(x)\big)) = \SP (\mathbb X(x)' \mathbb X^{-1}(x)) \\
= \SP (B(x)\sigma_2 B^*(x) \mathbb X^{-1}(x)) =
\tr (\sigma_2 B^*(x) \mathbb X^{-1}(x)B(x)).
\end{multline}
Since any two realizations of a symmetric function
are (weakly) isomorphic, one obtains from standard theorems \cite{bib:bgr} 
in realization theory of analytic at infinity functions that they will
have the same tau function up to a scalar. 
In other words, this notion is independent of the realization we choose for the given function $S(\lambda,x)$.

\subsection{\label{sec:Scatt}Scattering data versus $S(\lambda,0)$ and Gelfand-Levitan equation}
Following \cite{bib:FadeevInv} for the case $\int\limits_0^\infty x |q(x)| dx<\infty$, there are introduced Jost solutions
$\phi(x,s)$ and $f(x,s)$ \cite{bib:Jost, bib:Levinson}
\begin{eqnarray}
\label{eq:phiCond}	\phi(x,s) & : & \phi(0,s) =0, \quad \phi'(0,s) = 1, \\
\label{eq:fCond}		f(x,s) & : & \lim\limits_{x\rightarrow\infty} e^{-isx} f(x,s) = 1.
\end{eqnarray}
Defining further $M(k) = \phi'(x,k) f(x,s) - f'(x,s) \phi(x,k)$ one reconstructs the potential $q(x)$
using the Gelfand-Levitan equation \cite{bib:GL} (or alternatively the Marchenko equation \cite{bib:Marchenko}). 
There are two steps essential for
this construction: one considers the case when the spectrum of the operator $L$ is purely continuous
and the case when this spectrum additionally contains finite number of points.
For the purely continuous case, one proves that there is a solution $K(x,y) $ of
the Gelfand-Levitan equation
\cite[(8.5)]{bib:FadeevInv}
\begin{equation} \label{eq:GelfandLevitan}
K(x,y) + \Omega(x,y) + \int\limits_0^x K(x,t) \Omega(t,y) dt = 0, \quad x>y,
\end{equation}
where $\Omega(x,y)$ is uniquely defined from $M(k)$ by \cite[(8.4)]{bib:FadeevInv}
\[ \Omega(x,y) = 2/\pi \int_0^\infty \dfrac{\sin(kx)}{k}[\dfrac{1}{M(k)M(-k)}-1] \dfrac{\sin(ky)}{k} k^2 dk.
\]
The formula for the potential is \cite[(10.4)]{bib:FadeevInv}
$q(x) = 2 \dfrac{d}{dx} K(x,x)$. Then one make a modification, so that the discrete spectrum is taken into
account \cite[(8.14, 8.15)]{bib:FadeevInv}. We will present analogues of these formulas in our setting.
Suppose that we are given an SL vessel \eqref{eq:DefSLV}
\[ \mathfrak{V} = (A, B(x), \mathbb X(x); \sigma_1(x), 
\sigma_2(x), \gamma(x), \gamma_*(x);
\mathcal{H}, \mathcal{E};\mathrm I), 
\]
and let us fix an arbitrary $x_0\in\mathrm I$. Define
\begin{equation} \label{eq:OmeagaDef}
\Omega(x,y) = \bbmatrix{1&0} B^*(x) \mathbb X^{-1}(x_0) B(y) \bbmatrix{1\\0}
\end{equation}
and
\begin{equation} \label{eq:KDef}
K(x,y) = -\bbmatrix{1&0} B^*(x) \mathbb X^{-1}(x) B(y) \bbmatrix{1\\0}.
\end{equation}
Then the Gel'fand-Levitan equation \eqref{eq:GelfandLevitan} holds
\[ \begin{array}{llllll}
K(x,y) + \Omega(x,y) + \int\limits_{x_0}^x K(x,t) \Omega(t,y) dt = \\
= K(x,y) + \Omega(x,y) - \int\limits_{x_0}^x \bbmatrix{1&0} B^*(x) \mathbb X^{-1}(x) B(t) \bbmatrix{1\\0}  \bbmatrix{1&0} B^*(t) \mathbb X^{-1}(x_0) B(y) \bbmatrix{1\\0}dt  \\
= K(x,y) + \Omega(x,y) - \bbmatrix{1&0} B^*(x) \mathbb X^{-1}(x)  \int\limits_{x_0}^x B(t) \sigma_2 B^*(t) dt \mathbb X^{-1}(x_0) B(y) \bbmatrix{1\\0}  \\
= \text{ using vessel condition \eqref{eq:DB} }  \\
= K(x,y) + \Omega(x,y) - \bbmatrix{1&0} B^*(x) \mathbb X^{-1}(x) (\mathbb X(x) - \mathbb X(x_0)) \mathbb X^{-1}(x_0) B(y) \bbmatrix{1\\0}  \\
= K(x,y) + \Omega(x,y) - \bbmatrix{1&0} B^*(x) \mathbb X^{-1}(x_0) B(y) \bbmatrix{1\\0} +
\bbmatrix{1&0} B^*(x) \mathbb X^{-1}(x)  B(y) \bbmatrix{1\\0} = 0.
\end{array} \]
Finally, the formula \eqref{eq:qbeta} for the potential gives
\begin{equation} \label{eq:qFromK}
 \begin{array}{llllll}
q(x) = 2\dfrac{d}{dx} \SP ( \mathbb X^{-1}(x)  \dfrac{d}{dx} \mathbb X(x)) = 
\SP ( \mathbb X^{-1}(x)  B(x) \sigma_2 B^*(x)) = \\
\quad = 2\dfrac{d}{dx} \SP ( \mathbb X^{-1}(x)  B(x) \bbmatrix{1\\0} \bbmatrix{1&0} B^*(x)) = 
\dfrac{d}{dx}\big( \bbmatrix{1&0} B^*(x) \mathbb X^{-1}(x)  B(x) \bbmatrix{1\\0} \big) = \\
\quad = - 2\dfrac{d}{dx} K(x,x),
\end{array} \end{equation}
which is identical to \cite[(10.4)]{bib:FadeevInv}.

\section{KdV evolutionary vessels}
To give motivation for the definition of a KdV evolutionary vessel, we present a
realization of the 3-soliton \cite[Example 3.3]{bib:BKSWSolitons}. 
Choose $\mathcal H = \mathbb C^3$ and consider the collection, defining a vessel, where the
operators $B, \mathbb X$ depend on the variable $t$:
\[
\mathfrak{E}_{KdV} = (A, B(x,t), \mathbb X(x,t); \sigma_1, 
\sigma_2, \gamma, \gamma_*(x,t);
\mathbb C^3,\mathbb C^2).
\]
Suppose that the function $B(x,t)$ satisfies, in addition to \eqref{eq:DB}, 
the following differential equation
\begin{equation} \label{eq:DBt}
 \dfrac{\partial}{\partial t} B(x,t) = i A \dfrac{\partial}{\partial x} B(x,t)
\end{equation}
and that $\mathbb X(x,t)$ satisfies
\begin{equation} \label{eq:DXt}
\dfrac{\partial}{\partial t} \mathbb X = 
 i A B \sigma_2 B^* - i B\sigma_2 B^* A^* +
 i B\gamma B^*.
\end{equation}
Then for some nonzero constants $b_1,b_2,b_3$ and mutually different positive $k_i$'s we can solve all these
equations to obtain the following formulas:
\[ \begin{array}{lllll}
A & = \diag(-ik_j^2) = \bbmatrix{-ik_1^2&0&0\\0&-ik_2^2&0\\0&0&-ik_3^2} = -A^*, \\ 
B(x,t) & = \bbmatrix{e^{k_1x+k_1^3t} b_1&0&0\\0&e^{k_2x+k_2^3t} b_2&0\\0&0&e^{k_3x+k_3^3t} b_3} \bbmatrix{1&ik_1\\1&ik_2\\1&ik_3},\\
\mathbb X(x,t) & = I + [\dfrac{e^{(k_i+k_j)x+(k_i^3+k_j^3)t}}{k_i+k_j}b_ib^*_j] .
\end{array} \]
The main reason why this example is interesting is the following lemma.
\begin{lemma} The tau function of the vessel $\mathfrak{V}_{KdV}$ is
\begin{multline*}
\tau(x,t) = \det \mathbb X(x,t) = 1 + \sum\limits_{i=1}^3 c_i e^{2k_ix+2k_i^3t}  \\
+ \sum\limits_{1\leq i<j\leq 3} c_i c_j a_{ij} e^{2(k_i+k_j)x+2(k_i^3+k_j^3)t}  \\
+ c_1c_2c_3 a_{12}a_{13} a_{23} e^{2(k_1+k_2+k_3)x+2(k_1^3+k_2^3+k_3^3)t},
\end{multline*}
where $a_{ij}=\dfrac{(k_i-k_j)^2}{(k_i+k_j)^2}$ and $c_i=\dfrac{|b_i|^2}{2k_i}$.
\end{lemma}
\noindent\textbf{Proof:} Calculate the determinant, using the formula for Cauchy determinant.
\qed

Notice that the above formula coincides with the example of $3$-soliton in
\cite[Example 3.3]{bib:BKSWSolitons}.

\begin{defn} \label{def:KdvVessel}
A \textbf{KdV (evolutionary) vessel} is a collection of operators and spaces
\[
\mathfrak{V}_{KdV} = (A, B(x,t), \mathbb X(x,t); \sigma_1, 
\sigma_2, \gamma, \gamma_*(x);
\mathcal H,\mathcal{E}),
\]
where $\sigma_1, \sigma_2, \gamma$ are SL parameters (see Definition \ref{def:SLparam}). 
The operators are subject to
the vessel conditions \eqref{eq:DB}, \eqref{eq:XLyapunov}, \eqref{eq:DX}, 
\eqref{eq:Linkage}, \eqref{eq:DBt}, and \eqref{eq:DXt} and the normalization condition
\begin{equation} \label{eq:KdVNormalCond}
\tr(\sigma_1 B^*\mathbb X^{-1} B) = 0.
\end{equation}
\end{defn}
The normalization condition becomes much simpler in certain cases:
\[ \tr(\sigma_1 B^*\mathbb X^{-1} B) = 
\tr(B \sigma_1 B^*\mathbb X^{-1}) = \tr((A\mathbb X + \mathbb XA^*)X^{-1})
= \tr(A+A^*),
\]
if the last expression make sense.
Here is an example of such a vessel, that implements solitons.
\begin{thm} \label{thm:Solitons}
Let $\sigma_1, \sigma_2, \gamma$ be defined as in Definition \ref{def:SLparam}.
Let $k_i$ ($i=1,\ldots,n$) be mutually different positive real numbers, and $b_i\in\mathbb C$
($i=1,\ldots,n$) nonzero.
Define the collection
\[
\mathfrak{V}_{KdV} = (A, B(x,t), \mathbb X(x,t); \sigma_1, 
\sigma_2, \gamma, \gamma_*(x);
\mathbb C^n,\mathcal{E}),
\]
where
\begin{equation} \label{eq:DefOperSoliton}
\begin{array}{lllll}
A & = \diag(-ik_j^2) = -A^*, \\ 
B(x,t) & = \diag(e^{k_jx+k_j^3t} b_j) \bbmatrix{1&ik_1\\\vdots&\vdots\\1&ik_n}, \\
\mathbb X(x,t) & = I_n + [\dfrac{e^{(k_i+k_j)x+(k_i^3+k_j^3)t}}{k_i+k_j} b_ib_j^*] 
\end{array} \end{equation}
and $\gamma_*$ is defined using \eqref{eq:Linkage}. Then $\mathfrak{V}_{KdV}$ is a KdV evolutionary vessel. 
The tau function of this vessel
$\tau(x,t) = \det(\mathbb X(x,t))$ is the tau function of
an n-soliton of the KdV equation.
\end{thm}
\noindent\textbf{Proof:} Let us show that the equations of the vessel hold. Equation
\eqref{eq:DB} is derived from the following calculations:
\[ \begin{array}{llll}
\dfrac{\partial}{\partial x} B(x,t) = \diag(k_j) B(x,t), \\
- (A B(x,t) \sigma_2 + B(x,t) \gamma) \sigma_1^{-1} =
\diag(ik_j^2) B(x,t) \bbmatrix{0&1\\0&0} - B(x,t) \bbmatrix{0&0\\i&0} \\
\hspace{1cm} = \diag(ik_j^2) \diag(e^{k_jx+k_j^3t} b_j) [ \bbmatrix{1&ik_1\\\vdots&\vdots\\1&ik_n}\bbmatrix{0&1\\0&0} + \diag(e^{k_jx+k_j^3t} b_j) \bbmatrix{1&ik_1\\\vdots&\vdots\\1&ik_n} \bbmatrix{0&0\\-i&0} = \\
\hspace{1cm} = \diag(e^{k_jx+k_j^3t} b_j) [\bbmatrix{0&ik_1^2\\\vdots&\vdots\\0&ik_n^2} +
\bbmatrix{k_1&0\\\vdots&\vdots\\k_n&0}] = \diag(e^{k_jx+k_j^3t} b_j) [\bbmatrix{k_1&ik_1^2\\\vdots&\vdots\\k_n&ik_n^2} = \\
\hspace{1cm} = \diag(k_j) B(x,t).
\end{array} \]
The Lyapunov equation \eqref{eq:XLyapunov} is derived as follows, where the first equality comes from the fact that
$A=-A^*$:
\[  \begin{array}{llll}
A \mathbb X(x,t) +  \mathbb X(x,t) A^* + B(x,t) \sigma_1 B^*(x,t) = \\
= A (I + [\dfrac{e^{(k_i+k_j)x+(k_i^3+k_j^3)t}}{k_i+k_j}b_ib^*_j]) +
(I + [\dfrac{e^{(k_i+k_j)x+(k_i^3+k_j^3)t}}{k_i+k_j}b_ib^*_j]) (-A) +
B(x,t) \sigma_1 B^*(x,t)  \\
= \diag(-ik_j^2)[\dfrac{e^{(k_i+k_j)x+(k_i^3+k_j^3)t}}{k_i+k_j}b_ib^*_j] +
[\dfrac{e^{(k_i+k_j)x+(k_i^3+k_j^3)t}}{k_i+k_j}b_ib^*_j] \diag(ik_j^2)  \\
\hspace{2cm} + \diag(e^{k_jx+k_j^3t} b_j) \bbmatrix{1&ik_1\\\vdots&\vdots\\1&ik_n} \sigma_1
\bbmatrix{1&\cdots&1\\-ik_1&\cdots & -ik_n} \diag(e^{k_jx+k_j^3t} b_j)  \\
= [\dfrac{e^{(k_i+k_j)x+(k_i^3+k_j^3)t}}{k_i+k_j}b_ib^*_j (\sqrt{-1}(-k_i^2+k_j^2)] +
[e^{(k_i+k_j)x+(k_i^3+k_j^3)t}b_ib^*_j \sqrt{-1} (k_i-k_j)]  \\
= [e^{(k_i+k_j)x+(k_i^3+k_j^3)t}b_ib^*_j (\dfrac{\sqrt{-1}(-k_i^2+k_j^2)}{k_i+k_j} +
\sqrt{-1} (k_i-k_j))] = 0. 
\end{array} \]
Equations \eqref{eq:DX} and \eqref{eq:DBt} are immediate. Finally, the fact that
the tau function is the tau function of an $n$-soliton follows from the more general Theorem
\ref{thm:KdVequation}. The normalization condition \eqref{eq:KdVNormalCond} is also immediate in this case.
\qed
\begin{thm}\label{thm:KdVequation}
Let $\mathfrak{V}_{KdV}$ be a KdV vessel, then the potential $q(x,t)$
of the output SL equation \eqref{eq:OutCC} satisfies KdV equation \eqref{eq:KdV}:
\[ q'_t = - \dfrac{3}{2} q q'_x + \dfrac{1}{4} q'''_{xxx}.
\]
\end{thm}
\noindent\textbf{Proof:} By \eqref{eq:qbeta}, the KdV equation for $q(x,t)$ follows from the
differential equation for $\beta(t,x)$:
\begin{equation} \label{eq:PDEbeta}
 4 \beta'_t = -6 (\beta'_x)^2 + \beta'''_{xxx}.
\end{equation}
Starting from $\dfrac{\tau'}{\tau}=-\beta$ and \eqref{eq:GKform}, we obtain that
\[ \beta(t,x) = - \bbmatrix{1&0} B^*(x,t) \mathbb X^{-1}(x,t) B(x,t) \bbmatrix{1\\0}.
\]
The next formula follows from \eqref{eq:DB}, \eqref{eq:XLyapunov}, and \eqref{eq:DX}:
\[ \sigma_1 (B^* \mathbb X^{-1})'=\sigma_2 B^* \mathbb X^{-1} A + \gamma_* B^* \mathbb X^{-1},
\]
for $\gamma_*=\bbmatrix{-i(\beta'-\beta^2)&-\beta\\\beta&i}$, defined from \eqref{eq:Linkage}. Using this, we find that
\[ \beta'_x = \tr(\bbmatrix{0&-i\\i&0} B^*\mathbb X^{-1} B) + \beta^2.
\]
Then, using properties of the trace, it follows that
\[ \beta''_{xx} = \tr(\bbmatrix{-2i&0\\0&0}B^*\mathbb X^{-1}AB) + 
\tr(\bbmatrix{0&-2i\beta\\4i\beta&-2}B^*\mathbb X^{-1}B) + 2 \beta^3.
\]
Finally,
\[ \beta'''_{xxx} = \tr(\bbmatrix{-8i\beta&4\\-4&0}B^*\mathbb X^{-1}AB) + 
\tr(\bbmatrix{0&4i\beta'-4i\beta^2\\0&-4\beta}B^*\mathbb X^{-1}B) + 
6 (\beta')^2.
\]
On the other hand, from the vessel conditions it follows that
\[ \begin{array}{lll}
\beta'_t & = & \tr(\bbmatrix{-2i\beta&1\\-1&0} B^*\mathbb X^{-1}AB) 
-\beta\tr(\bbmatrix{0&0\\0&1}B^*\mathbb X^{-1}B) + 
i (\beta'-\beta^2)\tr(\bbmatrix{0&1\\0&0}B^*\mathbb X^{-1}B)  \\
&  & - i \beta \tr(\sigma_2B^*\mathbb X^{-1}B\sigma_1 B^*\mathbb X^{-1}B).
\end{array} \]
Notice that the last term in this expression is zero by the normalization condition \eqref{eq:KdVNormalCond}, since
\[ \begin{array}{lll}
\tr(\sigma_2B^*\mathbb X^{-1}B\sigma_1 B^*\mathbb X^{-1}B) & =
-\beta \tr(\sigma_1 B^*\mathbb X^{-1}B) = 0.
\end{array} \]
Comparing the last two expressions, we obtain the formula \eqref{eq:PDEbeta}. \qed

\subsection{The class of transfer functions}
Before we define the class of the transfer functions of evolutionary KdV vessels $\mathfrak V_{KdV}$,
we prove the following proposition:
\begin{prop} \label{prop:TranFunKdV}
Let $\mathfrak E$ be an evolutionary KdV vessel in Definition \ref{def:KdvVessel}. Then
its transfer function
\begin{equation} \label{eq:KdVS}
 S(\lambda,x,t) = I - B^*(x,t) \mathbb X^{-1}(x,t)(\lambda I - A)^{-1} B(x,t) \sigma_1
\end{equation}
has the following properties:
\begin{enumerate}
	\item $S(\lambda,x,t)$ is defined on the set 
	$\mathrm D_S = \mathrm J\times\mathrm I_x\times I_t$, 
	where $\mathrm I_x,\mathrm I_t$ are intervals and
	$\mathrm J=\SPEC(A)\subseteq\mathbb C$,
	\item For all $x\in\mathrm I_x, t\in\mathrm I_t$, 
	$S(\lambda, x,t)$ is an analytic function of $\lambda$ outside $\mathrm J$, and
	$S(\infty, x,t) = I$, when $\lambda$ approaches infinity, away from $\mathrm J$,
	\item For all $\lambda\not\in\mathrm J$, $S(\lambda, x,t)$ 
	is a differentiable function of $x$ and of $t$,
	\item $S(\lambda,x)$ satisfies the symmetry condition 
	\begin{equation} \label{eq:Symmetry}
		S^*(-\bar\lambda,x,t) \sigma_1 S(\lambda,x,t) =  \sigma_1
	\end{equation}
	for all $(\lambda,x,t)\in \mathrm D_S$,
	\item Multiplication by $S(\lambda,x,t)$ maps solutions $u(\lambda,x,t)$ 
	of the \textbf{input} LDE with the spectral parameter $\lambda$:
	\begin{equation} \label{eq:KdVInCC}
		-\sigma_1\dfrac{\partial}{\partial x}u(\lambda,x,t) + (\sigma_2 \lambda + \gamma(x,t))u(\lambda,x,t) = 0
	\end{equation}
	to solutions $y(\lambda,x,t)=S(\lambda,x,t)u(\lambda,x,t)$ of the 
	\textbf{output}	LDE with the same spectral parameter:
	\begin{equation} \label{eq:KdVOutCC}
		-\sigma_1(\dfrac{\partial}{\partial x} y(\lambda,x,t) + (\sigma_2 \lambda + \gamma_*(x,t))y(\lambda,x,t) = 0.
	\end{equation}
\end{enumerate}
\end{prop}
\noindent\textbf{Proof:} These properties are easily checked, and follow from formula
\eqref{eq:KdVS} for $S(\lambda,x,t)$. The first and second properties are immediate.
The third property follows from the differentiability assumptions
on the operators $\mathbb X(x,t), B(x,t)$. The fourth property follows from straightforward
calculations using the Lyapunov equation \eqref{eq:XLyapunov}:
\[ \begin{array}{lllll}
S^*(\mu, x,t) \sigma_1 S(\lambda, x,t) - \sigma_1 = \\
~~~~~~ - (\bar\mu+\lambda) \sigma_1
B^*(x,t) (\bar\mu I - A^*)^{-1}\mathbb X^{-1}(x,t) (\lambda I - A)^{-1}  B(x,t) \sigma_1 = 0
\end{array} \]
for $\mu=-\bar\lambda$. The fifth property follows
directly from the definitions (by substituting $y(\lambda,x,t) = S(\lambda, x,t) u(\lambda,x,t)$
in \eqref{eq:KdVOutCC} and using \eqref{eq:KdVInCC} for $u(\lambda,x,t)$), and from the formula \eqref{eq:KdVS} for $S(\lambda, x,t)$ (by using vessel conditions in order to differentiate it). \qed

\begin{defn} The class $\IKdV$ is defined as the class of functions 
$S(\lambda,x,t)$, defined on a set $\mathrm D_S = \mathrm J\times\mathrm I_x\times I_t$, 
where $\mathrm I_x,\mathrm I_t$ are intervals and
$\mathrm J\subseteq\mathbb C$ is a closed subset, possessing the properties of the Proposition
\ref{prop:TranFunKdV}. The function is called \textbf{regular}, if it is analytic
at infinity for all $x\in\mathrm I_x, t\in\mathrm I_t$.
\end{defn}
In the regular case (i.e. all operators are bounded and the transfer function
is analytic at infinity), similar to \cite[Theorem 2.7]{bib:GenVessel}, we obtain the 
following realization theorem.
\begin{thm}[Realization in regular case] \label{thm:KdVRealGen}
Given a transfer function $S(\lambda,x)\in\IKdV$, for which additionally the set 
$\mathrm J$ is bounded, there exists a regular KdV evolutionary vessel
$\mathfrak V_{KdV}$ whose transfer function coincides with $S(\lambda,x,t)$, probably on
smaller intervals $\mathrm I'_x\subseteq\mathrm I_x, \mathrm I'_t\subseteq\mathrm I_t$.
\end{thm}
\noindent\textbf{Proof:} Identical to \cite[Theorem 2.7]{bib:GenVessel}. \qed

Notice that, given a function $S(\lambda,x,t)\in\IKdV$, we can take a function
\begin{equation} \label{eq:HCommPhi}
H = \bbmatrix{h_1(\lambda) & \dfrac{i h_3(\lambda)}{\lambda} \\ h_3(\lambda) & h_1(\lambda)},
\end{equation}
which commutes with the fundamental matrix of the input LDE. As a result
$S(\lambda,x,t)H \in \IKdV$ and will correspond to the same vessel parameters (see
\cite{bib:SchurVessels} for the details). This turns out to be the only way of generating
a function, corresponding to the same vessel parameters.
\begin{thm}[Uniqueness in regular case]\label{thm:Uniquness}
Given KdV parameters $\sigma_1, \sigma_2, \gamma,\gamma_*(x,t)$, there exists a unique
initial value $S(\lambda,0,0)$ up to a symmetric, identity at
infinity matrix function $H(\lambda)$, defined in \eqref{eq:HCommPhi}.
\end{thm}
\noindent\textbf{Proof:} Identical to \cite[Lemma 4.3]{bib:SchurVessels}. \qed

\section{\label{sec:Examples}Examples of solutions of the KdV equation}
In this section we show different classes of solutions of the KdV equation \eqref{eq:KdV} arising
from different choices of evolutionary KdV vessels. In Table \ref{table:KdV} we presented
classes of solutions, arising from different types of evolutionary vessels. 
In the periodic and almost-periodic odd cases it is shown in the next sections that the correspondence is 1-1.
Adjusting the classical inverse scattering theory, we can show that if a solution of the KdV equation
satisfies at $t=0$ the condition $\int_{\mathbb R} (1+|x|) q(x) dx < \infty $ and is three times differentiable,
it can be realized by a KdV vessel.

The last pages are devoted to a proof of the correspondence between classes of solutions of
\eqref{eq:KdV}, appearing in column one and the corresponding evolutionary vessels, appearing
in column two. We consider only the case of Dissipative operators $\mathbb X_0=I$ and solutions 
of \eqref{eq:KdV}, which are odd in $x$. Moreover, we take the operator $A$ to be of a very
simple ``diagonal'' form. 

\subsection{Analytic--exponential class}
We first make a simple observation, related to the regular vessels.
\begin{thm} \label{thm:AnalyticParam}
Suppose that we are given a regular KdV vessel on an interval $\mathrm I$.
Then $\beta(x,t)$ can be extended to an analytic function in the neighborhood of $x_0$
on the interval $\mathrm I_x$ uniformly in $t\in\mathrm I_t$.
\end{thm}
\noindent\textbf{Proof:} Let us realize the transfer function $S(\lambda,t_2)$ for a fixed
value of $x_0, t_0$:
\[ S(\lambda,t_2^0) = I - B_0^* \mathbb X_0^{-1} (\lambda I - A)^{-1} B_0 \sigma_1,
\]
where all the operators are bounded. Let us define two fundamental matrices 
$\Phi(x,\lambda), \Psi(\lambda,t)$ by
\[\begin{array}{lll}
\dfrac{\partial}{\partial x} \Phi(x,\lambda) = -(\lambda \Phi(x,\lambda) \sigma_2+\Phi(x,\lambda)\gamma)\sigma_1^{-1},
\quad \Phi(0,\lambda) = I, \\
\dfrac{\partial}{\partial t} \Psi(t,\mu) = -i (\mu^2 \Psi(x,\mu) \sigma_2 + 
\mu \Psi(t,\mu)\gamma)\sigma_1^{-1}, \quad \Psi(0,\mu) = I.
\end{array} \]
Then due to \eqref{eq:DB} and \eqref{eq:DBt}, using the Danford-Schwartz calculus \cite{bib:DanSchw}
\[
B(x,t) = \oint\oint (\lambda I - A)^{-1} (\mu I -A)^{-1} B_0 \Phi(x,\lambda)\Psi(t,\mu)  d\lambda d\mu
\]
for $R > \|A\|$. Since the function $\Phi(x,\lambda), \Psi(t,\lambda)$ satisfies LDEs
with constant coefficients, it follows that $B(x,t)$ can be extended to an entire function in both variables $x,t$.
Similarly, one obtains that $\mathbb X(x,t)$ is an entire function of two variables.

Then $\mathbb X(x,t)$ is invertible in a neighborhood of $(0,0)$, since it
is invertible at $(x,t)=(0,0)$ and satisfies a differential equation. 
Finally, from the linkage condition \eqref{eq:Linkage} for $\gamma_*(x,t)$, it follows
that $\gamma_*(x,t)$ (and $\beta(x,t)$) are analytic in the same neighborhood.
\qed

Let us consider now an analytic potential $\beta(x,t)$. We will construct a transfer function
$S(\lambda,x,t)$ using its Taylor series in $\lambda$. The Taylor coefficients are called the \textit{Markov moments}
and we denote them by $H_n(x,t)$:
\begin{equation} \label{eq:SMoments}
 S(\lambda,x,t) = I -
\sum\limits_{n=0}^\infty \frac{1}{\lambda^{n+1}} H_n(x,t),
\quad H_n(x,t) = \bbmatrix{H_n^{11}(x,t) & H_n^{12}(x,t)\\H_n^{21}(x,t)&H_n^{22}(x,t)}.
\end{equation}
Then using the differential equation \eqref{eq:DS}, and comparing the coefficients of 
$\dfrac{1}{\lambda^n}$ we can obtain that (\cite[Section 5.6]{bib:SchurVessels})
\begin{equation} \label{eq:Hi+1HiFs}
\left\{ \begin{array}{llll}
H_{n+1}^{12} & = i H_n^{21} - \frac{\partial}{\partial x} H_n^{11} + \beta H_n^{11} , \\
H_{n+1}^{11}-H_{n+1}^{22} & = i( \frac{\partial}{\partial x} H_{n+1}^{12} - \beta H_{n+1}^{12}) , \\
\frac{\partial}{\partial x} (H_{n+1}^{11} + H_{n+1}^{22}) & = -  i (\beta'_x-\beta^2) H_{n+1}^{12} + \beta (H_{n+1}^{11} - H_{n+1}^{22}), \\
2i \frac{\partial}{\partial x} H_{n+1}^{21} & = \frac{\partial^2}{\partial x^2} H_{n+1}^{11} - 2 \beta \frac{\partial}{\partial x}H_{n+1}^{11}. 
\end{array} \right. \end{equation}
Rearranging the terms and performing simple algebraic manipulations, we can obtain that all the
entries depend on $H_n^{12}$ as follows:
\begin{eqnarray}
\label{eq:11From12} H_n^{11} & = \dfrac{i(H_n^{12})'}{2}  - i \int\limits_{0}^{x} \beta'(y) H_n^{12}(y) dy , \\
(H_n^{22})' & = \dfrac{i}{2} (H_n^{12})'' - i \beta (H_n^{12})' , \\
(H_n^{21})' & = \dfrac{1}{2} (H_n^{12})'''-\beta (H_n^{12})''-\beta' (H_n^{12})' + (2\beta\beta'-\beta'')H_n^{12}.
\end{eqnarray}
Differentiating the first equation in \eqref{eq:Hi+1HiFs} and substituting $H_n^{21}$ from the fourth row, 
we obtain that
\begin{equation} \label{eq:12From11}
 H_n^{12} = -\dfrac{(H_{n-1}^{11})'}{2} + \int\limits_0^x \beta'(y) H_{n-1}^{11}(y) dy.
\end{equation}
Define an operator $T$ on the space of analytic functions
\begin{equation} \label{eq:DefT}
T f(x) = 2 \int\limits_0^x \beta'(y) f(y) dy
\end{equation}
and define the semi-group, generated by $T$
\[ G'(x) = T G(x), \quad  G(0) = I,
\]
and denote $H_0^{12} = G K_0$, from \eqref{eq:11From12} it follows that $H_0^{11} = \dfrac{i}{2} G K_0'$.
Similarly, from \eqref{eq:12From11} it follows that $H_1^{12} = -\dfrac{i}{2^2} G K_0''$ and from \eqref{eq:11From12}
$H_1^{11}=-\dfrac{i^2}{2^3} G K_0'''$. By induction we can show that
\[ H_n^{12} = \dfrac{(-i)^n}{2^{2n}} G K_0^{(2n)}, \quad H_n^{12} = -\dfrac{(-i)^n}{2^{2n+1}} G K_0^{(2n+1)}.
\]
It is important that the moments will grow sub-exponentially uniformly in $x,t$ for the existence of the transfer
function $S(\lambda,x,t)$ as the following theorem states.
\begin{thm}\label{thm:AnalytBetaReal}
An analytic on $\mathrm I_x\times\mathrm I_t$ function $\beta(x,t)$ can be realized by a KdV vessel, if and only if
there exists $M, r >0$ such that on $\mathrm I_x\times\mathrm I_t$ for $K(x,t) = G^{-1} \beta(x,t)$ it holds that
\begin{equation} \label{eq:KBound} |\dfrac{\partial^n}{\partial x^n} K(x,t)| \leq M r^n.
\end{equation}
\end{thm}
\noindent\textbf{Proof:} The necessity of this condition follows from the preceding the theorem arguments,
since the moments are of the form $H_n(x,t) = B(x,t)\mathbb X^{-1} A^n B(x,t)$ and
\[ \| H_n \| \leq \| B(x,t)\|\|\mathbb X^{-1}\|\| A \|^n\|B(x,t) \| \leq M r^n,
\]
for $M' = \| B(x,t)\|\|\mathbb X^{-1}\|\|B(x,t) \|, r = \|A\|$. Therefore
\[ |K(x,t)| = \|G^{-1}\|H_n^{12}\| \leq \|G^{-1}\| M' r^n = M r^n,\quad M= \|G^{-1}\| M'.
\]
The converse also holds, since if the
condition \eqref{eq:KBound} holds one can define a function 
\[ 
S(\lambda,x,t) = \sum\limits_{n=0}^\infty \frac{1}{\lambda^{n+1}} H_n(x,t)
\]
which will be analytic at infinity $\lambda=\infty$, and satisfy the differential equation \eqref{eq:DS}.
In corollary \ref{cor:RegSol} we show that such a function can always be ``fixed'' to a symmetric function.
\qed

Let us concentrate now on the symmetry, which was not taken into account while constructing
the transfer function  in Theorem \ref{thm:AnalytBetaReal}.
Notice that the symmetry condition \eqref{eq:Symmetry}
can be checked at $x=0$ only, since the expression
\[ S^*(-\bar\lambda,x,t)\sigma_1S(\lambda,x,t)
\] 
has derivative zero with respect to $x$ and $t$. This follows immediately from the vessels
condition \eqref{eq:DB}, \eqref{eq:DX}, \eqref{eq:XLyapunov}, \eqref{eq:DBt}, and \eqref{eq:DXt}. So, it is enough to require
this condition at $(x,t)=(0,0)$. 
Taking a matrix of the form 
$H(\lambda) = \bbmatrix{h_1(\lambda) & \dfrac{i h_3(\lambda)}{\lambda} \\ h_3(\lambda) & h_1(\lambda)}$
\eqref{eq:HCommPhi}, we know from arguments preceeding Theorem \ref{thm:Uniquness} 
that $S(\lambda,x,t) H(\lambda)$ is in the class $\IKdV$ if and only if $S(\lambda,x,t) \in \IKdV$. Let us show that arbitrary matrix $K(\lambda)$ may be brought to 
a special form
which in turn may be brought to a symmetric one.
We use the notation $k^\bigstar$ \textit{for the reflection with respect to the imaginary axis}:
\[ k^\bigstar = k^*(-\bar\lambda).
\]
\begin{lemma} For an arbitrary matrix $K(\lambda) = \bbmatrix{k_1&k_2\\k_3&k_4}$, there exist the following 
decomposition
\[ \bbmatrix{k_1&k_2\\k_3&k_4} = \bbmatrix{s_1 & s_3-d\\s_3 & s_1} \bbmatrix{h_1&\dfrac{i}{\lambda}h_3\\h_3&h_1},
d=\dfrac{s_3s_1^\bigstar + s_1s_3^\bigstar}{s_1+s_1^\bigstar}.
\]
\end{lemma}
\textbf{Proof:} Choose the following parameters:
\[ \begin{array}{lll}
F = \dfrac{k_1-k_4}{k_2-\dfrac{i}{\lambda}k_3}, \quad
G = \dfrac{k_2-k_4F}{k_4-k_2\dfrac{i}{\lambda}F}, \\
H = -\dfrac{k_1}{k_4 F} (1+FG\dfrac{i}{\lambda})+G+\dfrac{1}{F}, \quad
T = \dfrac{k_1}{H^\bigstar(1+F(G-H))} \\
p^2 = \dfrac{TT^\bigstar}{T+T^\bigstar}(F+F^\bigstar) \dfrac{1}{HH^\bigstar} \quad (p=p^\bigstar),
\end{array} \]
Define next the entries as follows
\[
h_1 = \dfrac{T}{p}, \quad
h_3 = F h_1, \quad
s_1 = H^\bigstar p, \quad
s_3 = G s_1.
\]
Then by straightforward calculations we can show that the decomposition in the lemma holds.
Moreover, we can choose $k_1,k_4$ so that the behavior at infinity of all these functions is as follows
\[ F(\infty) = 0 \Rightarrow G(\infty)=0 \Rightarrow
H(\lambda) \approx - \dfrac{1}{F} \dfrac{-k_1+k_4}{K_4} = \dfrac{k_2-\dfrac{i}{\lambda} k_3}{k_4} \approx \dfrac{1}{\lambda}.
\]
Then
\[ T(\lambda) = \dfrac{k_1}{H^*(1+F(G-H)} \approx \lambda
\]
and
\[ p^2(\lambda) \approx \dfrac{(\lambda^2)}{2\lambda}(\dfrac{2}{\lambda}) \dfrac{1}{\dfrac{1}{\lambda^2}}= \lambda^2
\Rightarrow p(\lambda) \approx \lambda.
\]
As a result
\[ H p \approx \dfrac{1}{\lambda}\lambda = const.
\]
Consequently, $s_1$ is a constant at infinity and $s_3$ is zero.
\qed

\begin{lemma} \label{lemma:SSymmConstr}
Suppose that $K(\lambda)$ is a constant at $\lambda = \infty$ and is of the form
\[  \bbmatrix{s_1 & s_3-d\\s_3 & s_1},
d=\dfrac{s_3s_1^\bigstar + s_1s_3^\bigstar}{s_1+s_1^\bigstar}.
\]
Then defining $h_1(\lambda), h_3(\lambda)$ such that
\[ \left\{ \begin{array}{lll}
h_1 h_1^\bigstar + \dfrac{i}{\lambda} h_3 h_3^\bigstar = 
\dfrac{1}{1 + s_1 s_1^\bigstar}, \\
 h_1 h_3^\bigstar + h_3 h_1^\bigstar = 0,
\end{array} \right. \]
and the matrix $H = \bbmatrix{h_1(\lambda) & \dfrac{i h_3(\lambda)}{\lambda} \\ h_3(\lambda) & h_1(\lambda)}$,
we will obtain that $K(\lambda) H$ is symmetric.
\end{lemma} 
\noindent\textbf{Proof:} Immediate calculation, by substituting all the formulas into the symmetry condition
\[ K(\lambda) H(\lambda) \sigma_1  H^*(-\bar\lambda) K^*(-\bar\lambda) = \sigma_1.
\]
\qed

The condition \eqref{eq:KBound} can be translated on a condition on $\beta(x,t)$. Indeed, from the definition
of $K(x)$ it follows that $\beta(x,t) = G(x) K(x,t) = e^{T x} K(x,t)$, where the operator $T$ on the set of analytic
functions on $\mathrm I_x\times\mathrm I_t$ is defined by \eqref{eq:DefT}. 
Consequently, using the Leibniz rule for the differentiation we obtain
\begin{multline*}
|\beta^{(n)}(x,t)| \leq \sum\limits_{i=0}^n \binom{n}{i} \| T^i e^{Tx}\| |K^{(n-i)}(x,t)| \leq \\
\leq \sum\limits_{i=0}^n \binom{n}{i} \|T\|^i \|e^{Tx}\| Mr^{n-i} =
M(\|T\|+r)^n \|e^{Tx}\| \leq M L^n e^{b x} \leq M \dfrac{d^n}{dx^n} e^{k x},
\end{multline*}
for $L=\|T\|+r, b=\|T\|$ and $k=\max(L,b)$.
\begin{defn} The class of analytic functions $\beta(x,t)$ on $\mathrm I_x\times\mathrm I_t$, 
for which there exists constants $M,k>0$ such that
\begin{equation} |\beta^{(n)}(x,t)| \leq M \dfrac{d^n}{dx^n} e^{kx}
\end{equation}
is denoted by $\A(k,\mathrm I_x\times\mathrm I_t)$ and is called \textbf{analytic--exponential} class.
\end{defn}
\begin{cor}\label{cor:RegSol}
Every function $\beta(x,t)\in\A(k,\mathrm I_x\times\mathrm I_t)$ can be realized using
a KdV evolutionary vessel.
\end{cor}
\noindent\textbf{Proof:} Construct \eqref{eq:SMoments} the 
transfer function $S(\lambda,x,t)$ using the moments $H_n(x,t)$ satisfying the conditions
\eqref{eq:Hi+1HiFs}.
Then by Theorem \ref{thm:AnalytBetaReal} it satisfies the differential equation \eqref{eq:DS}
and using Lemma \ref{lemma:SSymmConstr} it can be brought to a symmetric form.
Thus $S(\lambda,x,t)$ belongs to the
class $\IKdV$ and is regular (i.e. analytic in $\lambda$ at infinity).
By Theorem \ref{thm:KdVRealGen} $S(\lambda,x,t)$ can be realized as the transfer function
of a KdV evolutionary vessel. \qed

\subsection{Continuous case: $\mathcal H = L^2(\Gamma)$}
In the classical scattering theory \cite{bib:FadeevInv},
under assumption $\int_0^\infty x |q(x)|dx<\infty$, it follows that the spectrum
of the corresponding operator $-\dfrac{d^2}{dx^2} + q(x)$ consists of a finite number of points and a cut along the imaginary axis. Since we are dealing with the pure continuous spectrum, we will refer the potentials corresponding to
$\Gamma\subseteq i\mathbb R^+$ as \textit{classical potentials}.

Let us first delineate the construction of a vessel on the curve $\Gamma=i\mathbb R_+$.
We define $\mathcal H = L^2(\Gamma)$ and the operators of the vessel as follows
(using notation $\mu=is^2, \delta = ir^2$):
\begin{equation} \label{eq:DefOperOddCont}
\begin{array}{lllll}
A & = \diag(i\mu) = -A^*, \\
B(x,t) & = b(\mu) \bbmatrix{\dfrac{\sin(sx-s^3t)}{s}& -i\cos(sx-s^3t)}, \\
(\mathbb X(x,t) f)(\mu) & = f(\mu) + \int\limits_0^\infty
\dfrac{\dfrac{\sin(sx-s^3t)}{s}\cos(r x-r^3t) - \cos(sx-s^3t)\dfrac{\sin(rx-r^3t)}{r}}{s^2-r^2} b(\mu)b(\delta)^* f(\delta) dr
\end{array} \end{equation}
It is important to notice that the involved operators are well-defined. The functions $\cos(sx-s^3t)$ and
$\dfrac{\sin(sx-s^3t)}{s}$ are actually functions of $s^2=-i\mu$. The operator $\mathbb X(x,t)$ acts on functions
$f(\mu)\in L^2(\Gamma)$ and gives as the output function $(\mathbb X(x,t)f)(\mu)$.

\begin{lemma} The collection $\mathfrak E$, defined by the operators in \eqref{eq:DefOperOddCont}
is a KdV evolutionary vessel.
\end{lemma}
\noindent\textbf{Proof:} Similar to the proof of Theorem \ref{thm:Solitons}.
\qed
\begin{lemma} The following equality holds for any $g(\mu)\in\mathcal H$
\[ \mathbb X(x,0) [g(\mu) \dfrac{\sin(sx)}{s}] = g(\mu) \dfrac{\sin(sx)}{s}.
\]
\end{lemma}
\noindent\textbf{Proof:} Notice that the right hand side of the equality satisfies the following
differential equation with initial conditions (recall that $\mu=is^2$):
\[ \dfrac{\partial^2}{\partial x^2} Y(\mu,x) + i\mu Y(\mu,x) = 0,
\quad Y(\mu,0) = 0, \quad Y'(\mu,0) = g(\mu).
\]
Let us show that the left hand side satisfies the same equation and the same initial conditions. Since $\mathbb X(0) = I$ and $\mathbb X'(0) = B(0) \sigma_2 B(0)^* = 0$,
we obtain that
\[ \mathbb X(0) [g(\mu) \dfrac{\sin(s0)}{s}] = 0, \quad
 \dfrac{\partial}{\partial x} \mathbb X(x) [g(\mu) \dfrac{\sin(sx)}{s}] |_{x=0}=
 \mathbb X'(0) [g(\mu) \dfrac{\sin(s0)}{s}] + \mathbb X(0) g(\mu) \cos(s0) = g(\mu), 
\]
which means that $\mathbb X(x) [g(\mu) \dfrac{\sin(sx)}{s}]$ satisfies the same initial conditions. Then
\[ \begin{array}{llllllll}
\dfrac{\partial^2}{\partial x^2} \mathbb X(x) [g(\mu) \dfrac{\sin(sx)}{s}] + i\mu \mathbb X(x) [g(\mu) \dfrac{\sin(sx)}{s}] = \\
= \mathbb X''(x) [g(\mu) \dfrac{\sin(sx)}{s}] + 2 \mathbb X'(x) [g(\mu) \cos(sx)] +
\mathbb X(x) [g(\mu) s \sin(sx)]  +
 i\mu \mathbb X(x) [g(\mu) \dfrac{\sin(sx)}{s}] \\
= \int\limits_0^L [B(\mu,x) \bbmatrix{0&i\\-i&0} B^*(\delta,x) g(\delta) \dfrac{\sin(tx)}{t} + 2 B(\mu,x) \sigma_2 B^*(\delta,x) g(\delta) \cos(tx) ]dt  \\
\hspace{1cm} + \int\limits_0^L [ \dfrac{B(\mu,x) \sigma_1 B^*(\delta,x)}{\mu-\delta} 
[ g(\delta) t \sin(tx)  + g(\delta) i \mu \dfrac{\sin(tx)}{t}] dt   \\
= \int\limits_0^L [B(\mu,x) \bbmatrix{0&i\\-i&0} B^*(\delta,x) g(\delta) \dfrac{\sin(tx)}{t} + 2 B(\mu,x) \sigma_2 B^*(\delta,x) g(\delta) \cos(tx) ]dt \\
\hspace{1cm} + \int\limits_0^L [ \dfrac{B(\mu,x) \sigma_1 B^*(\delta,x)}{\mu-\delta} 
g(\delta) i (\mu-\delta) \dfrac{\sin(tx)}{t}] dt  \\
= \int\limits_0^L [B(\mu,x) \bbmatrix{0&2i\\0&0} B^*(\delta,x) g(\delta) \dfrac{\sin(tx)}{t} + 2 B(\mu,x) \sigma_2 B^*(\delta,x) g(\delta) \cos(st) ]dt  \\
= \int\limits_0^L 2i c(\mu) \dfrac{\sin(sx)}{s} i \cos(tx) g(\delta) \dfrac{\sin(tx)}{t} +
2 c(\mu) \dfrac{\sin(sx)}{s} \dfrac{\sin(tx)}{t} g(\delta) \cos(tx) dt = 0. \qed
\end{array} \]
\begin{cor} The formula for the potential is
\begin{equation} \label{eq:qOdd}
 q(x) = 2 \int_0^L c(\mu) c^*(\mu) \dfrac{\sin(2sx)}{s} ds.
\end{equation}
\end{cor}
\noindent\textbf{Proof:} From the previous theorem, it follows that
\[ \mathbb X^{-1}(x) [c(\mu) \dfrac{\sin(sx)}{s}] = c(\mu) \dfrac{\sin(sx)}{s}.
\]
Then the result is immediate from formulas \eqref{eq:KDef} and \eqref{eq:qFromK}.
\qed
\begin{thm} There is a 1-1 correspondence between potentials $q(x)$ and
functions $c(\mu)$ satisfying \eqref{eq:qOdd}.
\end{thm}
From the uniqueness of the solution of the KdV equation we have the following theorem.
\begin{thm} \label{thm:L1Sol} There is 1-1 correspondence between $\beta(x,t)$ satisfying the
KdV equation \eqref{eq:KdV} and evolutionary KdV vessels, defined on the curve
$\Gamma=i\mathbb R_+$ using $c(\mu)$ defined from \eqref{eq:qOdd}.
\end{thm}

\subsection{Discrete case: $\mathcal H = \ell^2$}
In this section we present solutions for which the inner space $\mathcal H$ is discrete.
Suppose that the space $\mathcal H = \ell^2$ consists of 
infinite column-sequences ($(\cdot)^t$-denotes the transpose):
\[ \mathcal H = \ell^2 = 
\{ \bar{\mathbf x}_n^t = (x_1,x_2,\ldots)^t \mid \sum\limits_{n=1}^\infty |x_n|^2 < \infty \}.
\]
Let $\bar {\mathbf k}_n =(k_1,k_2,\ldots)$ be a sequence of mutually different numbers ($k_i\neq k_j$ for
$i\neq j$) and $\bar {\mathbf b}_n=(b_1,b_2,\ldots)\in\ell^2$ be another sequence. 
We define analogously to the finite dimensional case \eqref{eq:DefOperSoliton} the following operators:
\begin{equation} \label{eq:DefOperOddDiscrete}
\begin{array}{lllll}
A_p & = \diag(ik_n^2) = -A^*, \\ 
B_p(x,t) & = \diag(b_n) \bbmatrix{\dfrac{\sin(k_nx-k_n^3t)}{k_n}& -i\cos(k_nx-k_n^3t)}, \\
\mathbb X_p(x,t) & = I + [\dfrac{\dfrac{\sin(k_nx-k_n^3t)}{k_n}\cos(k_mx-k_m^3t) - \cos(k_nx-k_n^3t)\dfrac{\sin(k_mx-k_m^3t)}{k_m}}{k^2_n-k^2_m} b_nb_m^*] 
\end{array} \end{equation}
and the collection realizing an odd (for $t=0$) periodic solution of the KdV equation 
\[ \mathfrak E_{p,odd} = (A_p, B_p(x,t), \mathbb X_p(x,t); \sigma_1, 
\sigma_2, \gamma, \gamma_*(x,t);
\mathcal H,\mathbb C^2)
\]
\begin{lemma} The collection $\mathfrak E_{p,odd}$, defined by the operators in \eqref{eq:DefOperOddDiscrete}
is a KdV evolutionary vessel.
\end{lemma}
\noindent\textbf{Proof:} Similar to the proof of Theorem \ref{thm:Solitons}.
\qed

In order to reconstruct a vessel for a given potential, we first notice that for
$t=0$ it is similar to the classical inverse scattering problem. As in the continuous case $\mathcal H=L^2(\Gamma)$,
we can calculate the inverse $\mathbb X(x,0)$.
\begin{lemma} The following equality holds
\[ \mathbb X(x) [b_n \dfrac{\sin(k_nx)}{k_n}] = b_n \dfrac{\sin(k_nx)}{k_n}.
\]
\end{lemma}
\noindent\textbf{Proof:} Since $\mathbb X(0,t) = I$ and both functions satisfy the same differential equation,
the lemma follows.
\qed

As a result, we have the following corolalry:
\begin{cor} The formula for $\beta(x,0)$ is
\begin{equation} \label{eq:betaOdd}
 \beta(x,0) = \sum_n b_n b_n^* \dfrac{\sin^2(k_nx)}{k_n^2}.
\end{equation}
\end{cor}
Recall that in the KdV equation \eqref{eq:PDEbeta} one uses the third derivative of $\beta$. As a result,
we must demand that the third derivative of $\sum_n b_n b_n^* \dfrac{\sin^2(k_nx)}{k_n^2}$ exists.
Differentiating this expression three times, we obtain that the following sum must converge:
\[ \sum_n b_n b_n^* k_n \sin(2k_n x).
\]
Therefore the necessary and obviously sufficient condition for the existence of the vessel is that the sequence
$b_n \sqrt{k_n}$ is square summable.

Finally, we can use the uniqueness of the solution of the KdV equation in order to
obtain the following theorem.
\begin{thm}\label{thm:PeriodSol}
There is a 1-1 correspondence between KdV evolutionary vessels $\mathfrak{V}_{KdV,p}$ and solutions
$\beta(x,t)$ of \eqref{eq:PDEbeta}. The correspondence is implemented by the formula \eqref{eq:betaOdd}.
\end{thm}
\noindent\textbf{Proof:} For $t=0$, we will obtain that $\beta(x,0)$ is represented
by the formula \eqref{eq:betaOdd}. Let us construct a KdV evolutionary vessel $\mathfrak{V}_{KdV}$
using $\mathcal H=\ell^2$ and the sequence $\{b_n\}$, obtained in from the Fourier-sinus transform 
of $\beta(x,0)$ \eqref{eq:betaOdd}. Then the
$\beta_{KdV}(x,t)$ of the vessel $\mathfrak{V}_{KdV,p}$ satisfies the same KdV equation
\eqref{eq:KdV} as the given $\beta(x,t)$. By the uniqueness of the solution they coincide.
In other words, the vessel $\mathfrak{V}_{KdV}$ realizes $\beta(x,t)$.
\qed

Notice that following the same lines, we could obtain a KdV vessel, realizing an almost periodic, odd solution
of \eqref{eq:KdV}. For this we need to use the Bohr spectrum of the function $\beta(0,t)$. In this cacse the formula
\eqref{eq:betaOdd} becomes
\begin{equation} \label{eq:betaOddAP}
 \beta(x,0) = \sum_n b_n b_n^* \dfrac{\sin^2(r_nx)}{r_n^2}.
\end{equation}
for $b_n\sqrt{k_n}$ square summable, and $r_n\in\mathbb R$. A KdV vessel $\mathfrak{V}_{KdV,ap}$
realizing an almost periodic solution
is identical to the periodic case with $k_n$'s substituted by $r_n$'s. 
The analogue of Theorem \ref{thm:PeriodSol} is as follows.
\begin{thm}\label{thm:AlmostPeriodSol}
There is a 1-1 correspondence between KdV evolutionary vessels $\mathfrak{V}_{KdV,ap}$ and almost periodic odd solutions
$\beta(x,t)$ of \eqref{eq:PDEbeta}. The correspondence is implemented by the formula \eqref{eq:betaOddAP}.\end{thm}

\subsection{Explicit formula for $\beta(x,t)$ when $\mathcal H=\ell^2$}
In order to explicitly construct a solution of the KdV equation for $\beta(x,t)$
\eqref{eq:PDEbeta} we suppose that
\[ \beta(x,t) = \sum_n |b_n(t)|^2 \dfrac{\sin^2(k_nx-k_n^3 t)}{k_n^2}.
\]
Then substituting it into equation \eqref{eq:PDEbeta}, we will obtain the following equality
\[ 4 \sum_n \dfrac{|b_n(t)|^2}{dt} \dfrac{\sin^2(k_nx-k_n^3 t)}{k_n^2} =
-6 (\beta'_x(x,t))^2.
\]
In other words,
\begin{multline} \label{eq:dbnt}
 4 \sum_n \dfrac{d|b_n(t)|^2}{dt} \dfrac{\sin^2(k_nx-k_n^3 t)}{k_n^2}
= -6 \sum_n\sum_m \dfrac{|b_n(t)b_m(t)|^2}{k_nk_m} \sin(2k_nx-2k_n^3 t)\sin(2k_mx-2k_m^3 t) \\
= -6 \sum_n\sum_m \dfrac{|b_n(t)b_m(t)|^2}{k_nk_m} 
\dfrac{\cos(2(k_n-k_m)x - 2(k_n^3-k_m^3)t) - \cos(2(k_n+k_m)x - 2(k_n^3+k_m^3)t)}{2}.
\end{multline}
Notice that
\[ k_n^3\pm k_m^3 = (k_n\pm k_m)^3 \mp 3k_nk_m(k_n\pm k_m).
\]
Using the formula for $\cos(x+y)$ equation \eqref{eq:dbnt} becomes
\[ \begin{array}{lllll}
4 \sum_n \dfrac{d|b_n(t)|^2}{dt} \dfrac{\sin^2(k_nx-k_n^3 t)}{k_n^2} & =
-3 \sum_n\sum_m \dfrac{|b_n(t)b_m(t)|^2}{k_nk_m} \big[ \\
&\cos\big(2(k_n- k_m)x - 2(k_n - k_m)^3t\big)  \cos\big(6k_nk_m(k_n- k_m) t\big)  \\
&\quad - \sin\big(2(k_n- k_m)x - 2(k_n- k_m)^3t\big)  \sin\big(6k_nk_m(k_n- k_m) t\big)  \\
&\quad \quad -
\cos\big(2(k_n+ k_m)x - 2(k_n+ k_m)^3t\big)  \cos\big(6k_nk_m(k_n+ k_m) t\big)  \\
&\quad \quad\quad +
\sin\big(2(k_n+ k_m)x - 2(k_n+ k_m)^3t\big)  \sin\big(6k_nk_m(k_n+ k_m) t\big)\big].
\end{array} \]
Assume next that $-k_m = k_{-m}$ (i.e. that there is an antisymmetry 
condition on the numbers $k_m$). Then
\[ \begin{array}{lllll}
\sum_n\sum_m \dfrac{|b_n(t)b_m(t)|^2}{k_nk_m} \big[
\sin\big(2(k_n - k_m)x - 2(k_n - k_m)^3t\big)  \sin\big(6k_nk_m(k_n+ k_m) \big]  \\
= \sum_n\sum_{-m} \dfrac{|b_n(t)b_{-m}(t)|^2}{k_nk_{-m}} \big[
\sin\big(2(k_n - k_{-m})x - 2(k_n -{-k_m})^3t\big)  \sin\big(6k_nk_{-m}(k_n- k_{-m}) \big]  \\
= \sum_n\sum_m \dfrac{|b_n(t)b_{-m}(t)|^2}{-k_nk_m} \big[
\sin\big(2(k_n + k_m)x - 2(k_n + k_m)^3t\big)  \sin\big(-6k_nk_m(k_n+ k_m) \big]  \\
= \sum_n\sum_m \dfrac{|b_n(t)b_{-m}(t)|^2}{-k_nk_m} \big[
\sin\big(2(k_n + k_m)x - 2(k_n +k_m)^3t\big)  \sin\big(6k_nk_m(k_n+ k_m) \big].
\end{array} \]
Therefore, if $|b_{-m}|=|b_m|$, we obtain that the terms involving $\sin$ are canceled.
As a result we obtain in the last formula that
\begin{multline*}
4 \sum_n \dfrac{d|b_n(t)|^2}{dt} \dfrac{\sin^2(k_nx-k_n^3 t)}{k_n^2} \\ 
= -3 \sum_n\sum_m \dfrac{|b_n(t)b_m(t)|^2}{k_nk_m} 
\cos\big(2(k_n+ k_m)x - 2(k_n+ k_m)^3t\big)  \cos\big(6k_nk_m(k_n+ k_m) t\big).
\end{multline*}
Finally, it remains to demand that the set $\Gamma=\{k_n\}$ has the property
$\Gamma+\Gamma = \Gamma$. We shall obtain that
\begin{multline*}
4 \sum_N \dfrac{d|b_N(t)|^2}{dt} \dfrac{1- \cos(2k_Nx-2k_N^3 t)}{2 k_N^2}
= -3 \sum_{k_n+k_m=k_N} \dfrac{|b_n(t)b_m(t)|^2}{k_nk_m} \cos\big(6k_nk_mk_N t\big)
\cos\big(2k_Nx - 2k_N^3t\big).
\end{multline*}
In other words, $b_N(t)$ must satisfy the differential equation
\begin{equation}\label{eq:dbNt}
\dfrac{1}{k_N^2} \dfrac{d|b_N(t)|^2}{dt} = -\dfrac{3}{2} \sum_{k_n+k_m=k_N} \dfrac{|b_n(t)b_m(t)|^2}{k_nk_m} \cos\big(6k_nk_mk_N t\big)
\end{equation}
and the normalization condition
\begin{equation}\label{eq:sumbN}
\sum_N \dfrac{d|b_N(t)|^2}{dt} \dfrac{1}{k_N^2} = 0.
\end{equation}
Notice that the condition $b_{-N}(t)=b_N(t)$ holds in this case, since
\[ \begin{array}{lllll}
\dfrac{1}{k_{-N}^2} \dfrac{d|b_{-N}(t)|^2}{dt} & = -\dfrac{3}{2} \sum_{k_n+k_m=k_{-N}} \dfrac{|b_n(t)b_m(t)|^2}{k_nk_m} \cos\big(6k_nk_mk_{-N} t\big)  \\
& = -\dfrac{3}{2} \sum_{k_{-n}+k_{-m}=k_N} \dfrac{|b_{-n}(t)b_{-m}(t)|^2}{k_{-n}k_{-m}} \cos\big(6k_{-n}k_{-m}k_N t\big)  \\
&= -\dfrac{3}{2} \sum_{k_n+k_m=k_N} \dfrac{|b_n(t)b_m(t)|^2}{k_nk_m} \cos\big(6k_nk_mk_N t\big) \\
& = \dfrac{1}{k_N^2} \dfrac{d|b_N(t)|^2}{dt}.
\end{array} \]
The equality follows from the initial condition $b_N(0) = b_{-N}(0)$, which follows in turn
from the oddness of $\beta(x,0)$.
Thus we obtain the following theorem.
\begin{thm} \label{thm:BetaForm} Let $\beta(x,t)$ be a solution of the KdV equation, possessing the property that
$\beta(x,0)$ is an odd function, such that \eqref{eq:betaOdd} holds for a set
$\Gamma=\{k_n\}$. Assume that $-\Gamma=\Gamma$ (i.e. $-k_N=k_{-N}$) and $\Gamma + \Gamma =
\Gamma$. Then,
\[ \beta(x,t) = \sum_N |b_N(t)|^2 \dfrac{\sin^2(k_Nx-k_N^3 t)}{k_N^2}
\]
is a solution of the integrated KdV equation \eqref{eq:PDEbeta} if and only if
$b_N(t)$ satisfies \eqref{eq:dbNt} and the normalization condition \eqref{eq:sumbN}.
\end{thm}

\section{Conclusions and Remarks}
One can study a wide range of questions using evolutionary vessels and their
generalizations. We list several questions that we believe have a great potential to contribute 
to the study of PDEs, corresponding to integrable systems.
\begin{enumerate}
	\item Which PDEs are obtained from the SL equation by changing the evolutionary conditions
	\eqref{eq:DBt} and \eqref{eq:DXt}?
	\item Can we develop a similar theory starting from Non-Linear Shr\" odinger Equation (NLS) 
	\cite{bib:SchurVessels}? The formulas
	developed in this work will look different, but the ideas will remain the same. The NLS
	vessels parameters are presented in the next definition.
	\begin{defn} Non-Linear Shr\" odinger (NLS) Equation parameters are given by
		\[ \begin{array}{lll}
		 \sigma_1 = \bbmatrix{1&0\\0&1},\quad
		\sigma_2 = \dfrac{1}{2} \bbmatrix{1&0\\0&-1},
		\quad \gamma =\bbmatrix{0 & 0 \\0 & 0}, \\
		\gamma_*(x,t) = \bbmatrix{0&\beta(x,t) \\-\beta^*(x,t)& 0}
		\end{array} \]
	\end{defn}
	\item How are the Segal-Wilson \cite{bib:SegalWilson} loop groups solutions of the KdV equation
	obtained using vessels? What is the space $\mathcal H$ in this case?
	\item There exist Krichever \cite{bib:Krich77} solutions of the KdV. 
	How can they be implemented using vessels on curves \cite{bib:GenVessel}? 
	\item Can one realize the Zakarov-Shabat scheme by generalizing these types of vessels to 
	infinite number of variables?
\end{enumerate}

\bibliographystyle{alpha}
\bibliography{../../biblio}

\end{document}